\documentclass[12pt]{article}
\usepackage{textcomp}
\usepackage{epsfig, color}
\usepackage{amssymb,amsmath,latexsym}
\usepackage{epstopdf}
\usepackage{mathrsfs}

\newtheorem{thm}{Theorem} %
\newtheorem{lem}{Lemma}
\newtheorem{conj}{Conjecture}
\newtheorem{prop}{Proposition}
\newtheorem{prob}{Problem}
\newtheorem{cor}{Corollary}

\newtheorem{corollary}[thm]{\rm\bfseries Corollary}
\newtheorem{claim}[thm]{Claim}
\newtheorem{proposition}[thm]{Proposition}

\newcommand {\red} {\textcolor{red}}

\newcommand \ARRAY[1]
{
\left \{
\begin{array}{ll}
#1
\end{array}
\right.
}

\newcommand{\proof}{{\noindent {\em Proof}.\quad}
\setcounter{countclaim}{0}\setcounter{countcase}{0}}

\newcommand{\proofn}{{\noindent {\em Proof}.\quad}}

\newcommand{\proofend}{{\hfill$\Box$}}

\newcounter{countcase}

\newcounter{countclaim}
\def\inclaim{\addtocounter{countclaim}{1}
{\noindent {\bf Claim \thecountclaim}: }}

\newcommand {\rebibitem}[1] {\bibitem{#1} \red{[*: #1]}}%draft

\def\rebibitem {\bibitem}  %take it to remove labels in references.

\begin{document}
%\linenumbers

\title{On the sizes of bipartite 1-planar graphs
\thanks{This work is supported by MOE-LCSM, School of
   Mathematics and Statistics, Hunan Normal University, Changsha, Hunan 410081,
   P. R. China }}
\author{ Yuanqiu Huang\thanks{Corresponding author} \\
{\footnotesize  Department of Mathematics, Normal University of
Hunan, Changsha 410081, P.R.China}\\
{\footnotesize hyqq@hunnu.edu.cn}\\
Zhangdong Ouyang\\
{\footnotesize Department of Mathematics, Hunan First Normal University , Changsha 410205, P.R.China} \\
{\footnotesize oymath@163.com }\\
Fengming Dong\\
{\footnotesize National Institute of Education, Nayang Technological University, Singapore} \\
{\footnotesize donggraph @163.com }\\
}

\date{}
\maketitle

\begin{abstract}
A graph is called $1$-planar if it admits a drawing in the plane such that each edge is crossed at most once.
Let $G$ be a bipartite 1-planar graph with $n$ ($\ge 4$) vertices
and $m$ edges.
Karpov showed that $m\le 3n-8$ holds for even $n\ge 8$ and
$m\le 3n-9$ holds for odd $n\ge 7$.
Czap, Przybylo and \u{S}krabul\'{a}kov\'{a}
proved that if the partite sets of $G$ are of sizes $x$ and $y$,
then $m\le 2n+6x-12$ holds for $2\leq x\leq y$,
and conjectured that $m\le 2n+4x-12$ holds for $x\ge 3$ and $y\ge 6x-12$.
In this paper, we settle their conjecture and our result
is even under a weaker condition $2\le x\le y$.
\iffalse
The authors in this paper also constructed a bipartite 1-planar graph $G$ with partite sets of sizes $x$ and $y$ with $x\geq 3$ and $y\geq 6x-12$ such that $|E(G)|\geq 2|V(G)|+4x-12$, and  posed the following conjecture: {\it for any integers $x$,$y$ such that $x\geq 3$ and $y\geq 6x-12$, every bipartite 1-graph $G$ with partite sets of sizes $x$ and $y$ has at most $2|V(G)|+4x-12$}.  Another result, due to D.V.Karpov [J.Math. Sci. 196 (2014), 737-746], is that   $|E(G)|\leq  3|V(G)|-8$ for any bipartite 1-planar graph $G$ with $|V(G)|\geq 4$.  In this paper, we prove that if $G$ is a bipartite 1-graph of  partite sets of sizes $x$ and $y$  with
$2\leq x \leq y$, then $E(G)|\leq 2|V(G)|+4x-12$. Therefore, our result greatly improves J.Czap et al's the upper bound ``$2|V(G)|+6x-12$'' to ``$2|V(G)|+4x-12$", and implies that the above conjecture is correct  and  the restriction ``$y\geq 6x-12$"  is redundant (it only needs that $2\leq x\leq y$). Our result on the upper bound is  better than D. V.Karpov's above when  $x\leq \frac{1}{3}(y+4)$. We also completely solve a problem  raised by Sopena
that how many edges must be removed  from a complete bipartite graph so as to  possibly produce a 1-planar graph.\fi
\end{abstract}

\noindent
%\textbf{AMS classification}: 05C07, 05C15, 05C50\\
{\bf Keywords}: bipartite graph,  drawing, face,  1-planar graph.

\maketitle

\iffalse
%%%%%% THIS PART MUST BE PLACED IMMEDIATELY AFTER THE \maketitle COMMAND
%%%%%% BACK TO ORIGINAL FOOTNOTES
\makeatletter
\renewcommand\@makefnmark%
{\mbox{\textsuperscript{\normalfont\@thefnmark)}}}
\makeatother
%%%%%%
\fi

\section{Introduction}

 A {\em drawing} of a graph
$G=(V,E)$ is a mapping $D$ that assigns to each vertex in $V$ a
distinct point in the plane and to each edge $uv$ in $E$ a
continuous arc connecting $D(u)$ and $D(v)$. We often make no
distinction between a graph-theoretical object (such as a vertex, or
an edge) and its drawing. All drawings considered here are
such ones  that no edge crosses itself, no
two edges cross more than once, and no two edges incident with the
same vertex cross.  The {\em crossing
    number} of a graph $G$ is the smallest number of
crossings in any drawing of $G$.

A drawing of a graph is 1-{\it planar} if each of its edges is
crossed at most once. If a graph has a 1-planar drawing, then it is
1-{\it planar}. The notion of 1-planarity was introduced in 1965 by
Ringel \cite{GR}, and since then many properties of 1-planar graphs have
been studied (e.g. see the survey paper \cite{SK}).

It is well-known that any simple planar graph with $n$ ($n\geq 3$) vertices
has at most $3n-6$ edges, and
a simple and bipartite graph  with $n$ ($n\geq 3$) vertices has at most $2n-4$ edges.
Determining the maximum number of edges in 1-planar graphs
with a fixed number of vertices
has aroused great interest of many authors (see, for example, \cite{HR}, \cite{JDT},\cite{IT}, \cite{JP}, \cite{XZ}).
It is known that \cite{HR,IT, JP} any 1-planar graph with $n$ $(\geq 3)$ vertices has at most $4n-8$ edges.
For bipartite 1-planar graphs,
an analogous result was due to Karpov \cite{DVK}.

\iffalse
\vskip 0.2cm
\noindent
{\bf Theorem 1} \cite{DVK}.
{\it Let $G$ be a bipartite 1-planar graph with $n$ vertices. Then $G$ has at most $3n-8$ edges for even $n\not=6$, and at most $3n-9$ edges for odd $n$ and for $n=6$. \red{For all $n\geq 4$, these bounds are tight}.
}
\fi

\begin{thm}[\cite{DVK}]\label{th1}
Let $G$ be a bipartite $1$-planar graph with $n$ vertices. Then $G$ has at most $3n-8$ edges for even $n\not=6$, and at most $3n-9$ edges for odd $n$ and for $n=6$.
For all $n\geq 4$, these bounds are tight.
\end{thm}

%\vskip 0.2cm
Note that Karpov's upper bound on the size of
a bipartite 1-planar graph is in terms of its vertex number.
When the sizes of partite sets
in a bipartite 1-planar graph are taken into account,
Czap, Przybylo and \u{S}krabul\'{a}kov\'{a}~\cite{CPS} obtained
another upper bound for its size
(i.e., Corollary 2 in \cite{CPS}).

\iffalse
\vskip 0.2cm
\noindent
{\bf Theorem 2 } \cite{CPS}. {\it If $G$ is a bipartite 1-planar graph such that the partite sets of $G$ have sizes $x$ and $y$, $2\leq x\leq y$, then $|E(G)|\leq 2|V(G)|+6x-16$. }
\fi

For any graph $G$, let $V(G)$ and $E(G)$ denote its vertex set
and edge set.

\begin{thm}[\cite{CPS}]\label{th2}
If $G$ is a bipartite $1$-planar graph
with partite sets of sizes $x$ and $y$,
where $2\leq x\leq y$, then $|E(G)|\leq 2|V(G)|+6x-16$.
\end{thm}

%\vskip 0.2cm

For each pair of integers $x$ and $y$
with $x\geq 3$ and $y\geq 6x-12$,
the authors in \cite{CPS} constructed a bipartite 1-planar graph $G$
with partite sets of sizes  $x$ and $y$
such that $|E(G)|=2|V(G)|+4x-12$ holds.
\iffalse
also investigated lower bounds on the number of edges of bipartite 1-planar graph with some restrictions on sizes of the partite sets;  they  \cite{CPS} constructed infinitely many bipartite 1-planar graphs $G$ with  partite sets of sizes  $x$ and $y$ , where $x\geq 3$ and $y\geq 6x-12$, such that $|E(G)|\geq 2|V(G)|+4x-12$,
\fi
Moreover,  they believed this lower bound is optimal
for such graphs and thus posed the following conjecture.

\iffalse
\vskip 0.2cm
\noindent
{\bf  Conjecture 3} \cite{CPS}. {\it For any integers $x$ and $y$ such that $x\geq 3$ and $y\geq 6x-12$, every bipartite 1-planar graph $G$ with partite sets of sizes $x$ and $y$ has at most $2|V(G)|+4x-12 $.}
\vskip 0.2cm
\fi

\begin{conj}[\cite{CPS}]\label{con1}
For any integers $x$ and $y$ with $x\geq 3$ and $y\geq 6x-12$,
if $G$ is a bipartite 1-planar graph with partite sets of sizes $x$ and $y$,
then $|E(G)|\le 2|V(G)|+4x-12$.
\end{conj}

In this paper we obtain the following result which proves
Conjecture~\ref{con1}.

\iffalse
\vskip 0.2cm
\noindent
{\bf Theorem 4}. {\it  Let $G$ be  a bipartite 1-graph  that has   partite sets of sizes $x$ and $y$   with
$2\leq x \leq y$. Then we have that  $|E(G)|\leq 2|V(G)|+4x-12$, and the upper bound is best possible.}
\fi

\begin{thm}\label{th-main}
Let $G$ be a bipartite $1$-planar graph with partite sets
of sizes $x$ and $y$, where  $2\leq x \leq y$.
Then $|E(G)|\leq 2|V(G)|+4x-12$,
and the upper bound is best possible.
\end{thm}

%\vskip 0.2cm

The result in \cite[Lemma 4]{CPS} shows that the upper bound for $|E(G)|$
in Theorem~\ref{th-main} is tight.
%The optimality on the upper bound in Theorem~\ref{th-main}  can be directly seen from  the constructions described  in \cite{CPS}.
%Theorem~\ref{th-main} shows that the conclusion of Conjecture~\ref{con1} holds with a more relaxed condition on $x$ and $y$: $2\leq x\leq y$.
%greatly improves J.Czap et al's  upper bound ``$2|V(G)|+6x-12$'' to ``$2|V(G)|+4x-12$", and  shows that the above conjecture is correct  but the restriction ``$x\geq 3$ and $y\geq 6x-12$"   is redundant (it only needs that $2\leq x\leq y$).
Also,  if $x\leq \frac{1}{3}(y+4)$,
Theorem~\ref{th-main} provides a better upper bound for $|E(G)|$
than Theorem~\ref{th1}.

The authors in \cite{CPS} mentioned  a question of Sopena \cite{S}:
{\it How many edges we have to remove from the complete bipartite graph with given sizes of the partite sets to obtain a 1-planar graph?}
It is not hard to see that Theorem~\ref{th-main}  implies
the follow corollary which answers the problem.
%has completely solved this  problem by  the  following .

\iffalse
\vskip 0.2cm
\noindent
{\bf Corollary 5}. {\it Let $K_{x, y}$ be a complete bipartite graph with sizes of the partite sets  $2\leq x\leq y$.  Then at least $(x-2)(y-6)$ edges  must be removed from $K_{x,y}$ such that the resulting graph becomes possibly a 1-planar graph, and the lower bound on the number of removed edges is best possible.}
\fi

\begin{cor}\label{cor1}
Let $K_{x, y}$ be the complete bipartite graph
with partite sets of sizes $x$ and $y$, where $2\leq x\leq y$.
Then at least $(x-2)(y-6)$ edges  must be removed from $K_{x,y}$ such that the resulting graph becomes possibly a 1-planar graph, and the lower bound on the number of removed edges is best possible.
\end{cor}

\iffalse
\vskip 0.2cm
\noindent
\textcolor{red}{{\bf Corollary 6}. {\it Let $G$  be a tripartite 1-planar graph $K_{x, y, z}$ with sizes of the partite sets  $2\leq x\leq y\leq z$.  Then $|E(G)|\leq 2|V(G)|+4x+6y+2z-24$.}}

 \vskip 0.2cm
 \fi

The remainder of this paper is arranged as follows.
In Section~\ref{sec2}, we explain some terminology and notation used in this paper.
In Section~\ref{sec3}, under some restrictions,
we present several structural properties on an extension of $D^{\times}$ for a 1-planar drawing $D$
of a bipartite 1-planar graph $G$,
where $D^{\times}$ is a plane graph introduced in Section~\ref{sec2}.
Some important lemmas for proving  Theorem~\ref{th-main}
are given in Section~\ref{sec4}, while
the proof of this theorem is completed in Section~\ref{sec5}.
Finally, we give some further problems in Section~\ref{sec6}.

\section{Terminology and notation\label{sec2}}

All graphs considered here are simple, finite and undirected,
unless otherwise stated.
For terminology and notation not defined here,
we refer to \cite{JAB}.
%We use $V(G)$ and $E(G)$ to denote the vertex set and the edge set of a graph $G$.
For any graph $G$ and $A\subseteq V(G)$,
let $G[A]$ denote the subgraph of $G$
with vertex set $A$ and edge set
$\{e\in E(G): e \mbox{ joins two vertices in }A\}$.
$G[A]$ is called the subgraph of $G$ induced by $A$.
For a proper subset $A$ of $V(G)$,
let $G-A$ denote the subgraph $G[V(G)\setminus A]$.
%the {\it vertex-induced subgraph} of a graph $G$ induced by $A$ is such one that is obtained from $G$ by deleting the vertices in $A$ together with their incident edges.

A {\it walk} in a graph is alternately a vertex-edge sequence; the walk is {\it closed } if its original vertex and terminal vertex are the same.  A {\it  path} (respectively, a {\it cycle }) of a graph is a walk (respectively, a closed walk) in which   all vertices  are distinct; the {\it length} of a path or cycle is the number of edges contained  in it. A path (respectively, a cycle) of length $k$  is  said to be a  {\it $k$-path} (respectively, {\it $k$-cycle}). If a cycle $C$  is  composed of two paths $P_1$ and  $P_2$,   we  sometimes write $C=P_1\cup P_2$.

 A {\it  plane graph} is a planar graph together with
a drawing without crossings,
and at this time  we  say that $G$ is embedded in the plane.
A plane graph $G$ partitions the plane into  a number of connected regions,  each of which is called  a {\it face }of $G$.  If a face is homeomorphic to an open disc, then it is called {\it cellular}; otherwise, {\it noncellular}.   Actually, a noncellular face is homeomorphic to an open disc with a few  removed  ``holes".  For a cellular face $f$,  the boundary of $f$  can be regarded as a closed walk of $G$, while for a noncellular  face $f$, its boundary  consists of  many  disjoint closed walks of $G$. The {\it  size} of a face is the number of the edges contained
in its the boundary with each repeated edge counts twice.
A face of size $k$ is also said to be a $k$-{\it face}.

It is known that a plane graph $G$ has no  noncellular faces if and only if $G$ is connected.
For a connected plane graph $G$, the well-known Euler's formula states that $|V(G)|-|E(G)|+|F(G)|=2$, where $F(G)$ denotes the face set of $G$.

A cycle $C$  of a plane graph $G$ partitions the plane into two open regions,  the bounded one (i.e., the {\it interior } of $C$)
and the  unbounded one (i.e., the {\it exterior} of $C$).
We denote  by $int(C)$ and  $ext(C)$  the interior and exterior of $C$, respectively, and their closures by $INT(C)$ and $EXT (C)$.
Clearly, $INT(C)\cap EXT(C) =C$.
A cycle $C$ of   a plane graph $G$ is said to be  {\it separating } if both $int(C)$ and $ext(C)$ contain  at least one  vertex of $G$.

Let $D$ be a 1-planar drawing of a graph $G$. The {\it associated plane graph} $D^{\times}$  is the plane graph that is obtained from $D$ by turning all crossings
of $D$ into new vertices of degree four;  these new vertices of degree four are called the {\it crossing vertices} of $D^{\times}$.

\vskip 0.3cm

\section{An extension of  $D^{\times}$ \label{sec3}
%on a bipartite 1-planar graph
}

%\vskip 0.2cm\noindent

Throughout this section, we always assume
that the considered graph $G$ (possibly disconnected) is a bipartite 1-planar graph with partite sets $X$ and $Y$, where $3\leq |X|\leq |Y|$.
Let $D$ be a 1-planar drawing of $G$ with the minimum number of crossings,  and $D^{\times}$ be the associated graph of $D$ with the crossing vertex set $W$.

\begin{figure}[!ht]
%\begin{center}
\centering
\includegraphics[width=10cm]{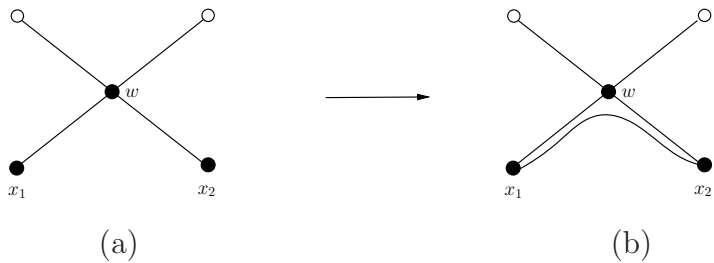} %{1-2}
%\vskip 0.2cm

(a) \hspace{6 cm} (b)

\caption{The extension of $D^{\times}$. \label{f1}}

%\centerline{\it Figure 1:  The extension of $D^{\times}$. }
%\end{center}
\end{figure}

Note that subsets $X,Y$ and $W$ form a partition of $V(D^{\times})$.
We color the vertices in $X,Y$ and $W$ by the black color,
white color and red color respectively.
As stated in \cite{CPS},  $D^{\times}$ can be extended to
a plane graph, denoted by $D^{\times}_{W}$,
by adding edges joining black vertices as described below:
\begin{quote}
for each vertex $w$ in $W$,
it is adjacent to two black vertices in $D^{\times}$, say $x_1$ and $x_2$,
and we draw an edge, denoted by $e_w$, joining $x_1$ and $x_2$
which is ``most near" one side of the path $x_1wx_2$ of $D^{\times}$
such that it does not cross with any other edge,
as shown in Figure~\ref{f1} (b).
\end{quote}

\iffalse
we extend $D^{\times} $ in the following way. Let $w$ be a crossing vertex in $W$ incident with two  black vertices $x_1$ and $x_2$. We draw a new  noncrossing edge $e$ by joining $x_1$ and $x_2$  ``most near"  one side of the path $x_1wx_2$ (see Figure 1).
\fi

\begin{figure}[!ht]

\centering
\includegraphics[width=15cm]{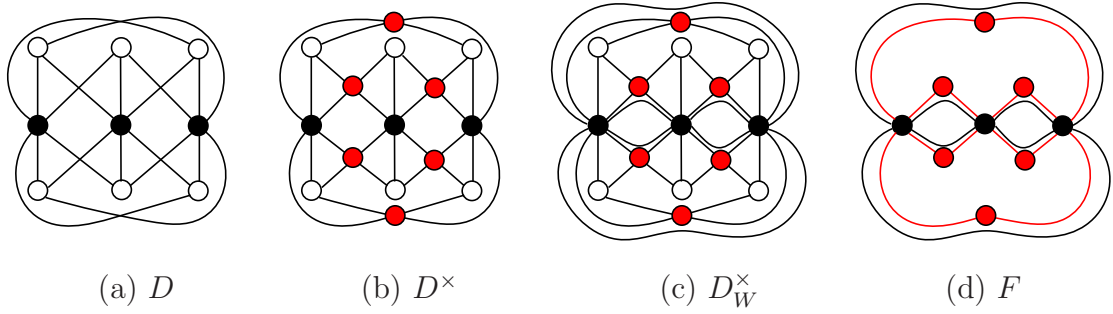}

(a) $D$ \hspace{2.2 cm} (b) $D^{\times}$ \hspace{2.4 cm}
(c) $D^{\times}_W$ \hspace{2.2 cm} (d) $F$

\caption{$D, D^{\times}$, $D^{\times}_W$ and $F$,
where $D$ is a 1-planar drawing of $K_{3,6}$ \label{f6}}
\end{figure}

Observe that
$D^{\times}_W$ is a plane graph with $D^{\times}$
as its spanning subgraph and
the edge set of $D^{\times}_W$ is the union of
$E(D^{\times})$ and $\{e_w: w\in W\}$.
Although $D^{\times}$  is a simple graph,
$D^{\times}_W$ might contain parallel edges
(i.e., edges with the same pair of ends),
as there may exist two edges in $\{e_w: w\in W\}$
with the same pair of ends.
An example is shown in Figure~\ref{f6} (c),
where $D$ is a $1$-planar drawing of $K_{3,6}$.

Let $F_D$ (or simply $F$) and
$H_D$ (or simply $H$) denote the subgraphs $D^{\times}_{W}[W\cup X]$
and $D^{\times}_{W}[X]$ respectively.
Obviously, $H$ is a  subgraph of $F$
and its edge set is $\{e_w: w\in W\}$,
while the edge set of $F$ is the union of
$E(H)$ and $\{wx_1,wx_2\in E(D^{\times}):
w\in W\ \&\ x_1,x_2\in X\}$.

All vertices in $H$ are black
and the edges in $H$ are also  called  {\it black edges}.
%Note that the edge set of $F$ is the union of $\{wx\in E(D^{\times}): w\in W\ \&\ x\in X\}$and $E(H)$ (i.e., $\{e_w: w\in W\}$).
%Thus, $|E(F)|=3|E(H)|$.
Clearly, $W$ is an independent set in $F$
and each vertex in $W$ (i.e., a red vertex) is of degree $2$ in $F$.
The edges in $F$ incident with red vertices are called
{\it red edges}.
Thus, each edge in $F$ is either black or red,
as shown in Figure~\ref{f6} (d).

We have  the following facts on $D^{\times}$, $F$ and $H$:
\begin{enumerate}
\item [(1)] $D^{\times}$, $F$ and $H$ are simultaneously embedded in the plane;
\item [(2)] $F$ and $H$ are obviously   loopless, but they are  possibly disconnected;
\item [(3)] $w\rightarrow e_w$ is
a bijection from $W$ to $E(H)$,
where $w$ is a red vertex,
and thus  the number of crossings of $D$ equals to $|E(H)|$;
and

%joining black vertices $x_1$ and $x_2$

% where $w$ is a red vertex;

\item [(4)] $e_w\rightarrow x_1wx_2$
is a bijection from $E(H)$ to the set of $2$-paths in $F$
whose ends are black,
where  $w$ is a red vertex and
$x_1$ and $x_2$ are the black vertices in $D^{\times}$
adjacent to $w$.
\end{enumerate}

Moreover we have the following propositions.

%\noindent

\begin{figure}[!ht]
\begin{center}
\includegraphics[width=10cm]{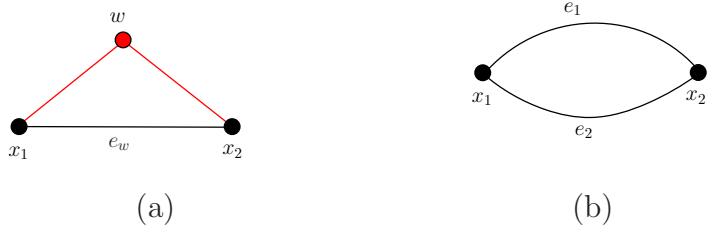}%{2-1}

(a) \hspace{5 cm} (b)

\caption{Some 3-cycles and 2-cycles in $F$.
\label{f3}}
\end{center}
\end{figure}

\begin{prop}\label{empty}
Let $e_w$ be an edge of $H$ with ends $x_1$ and $x_2$
and $C$ be the $3$-cycle of $F$  consisting  of  $e_w$  and  its corresponding  2-path $P=x_1wx_2$,
where $w$ is a red vertex
(see Figure~\ref{f3} (a)).
Then  $int (C)$   contains  none of  black vertices,  red vertices
and black edges in $F$;
in this sense we also say that $int(C)$ is ``empty".
\end{prop}

\proof By the definition of $D^{\times}_W$,
the proposition follows directly from the fact
that the drawing of edge $e_w$
is most near one side of the  2-path $x_1wx_2$ in $D^{\times}$
without crossings with edges in $D^{\times}$.
\proofend

\vskip 0.2cm

\begin{prop}\label{2-cycle}
Assume that $H$ has no separating 2-cycles.
If $C$ is  a 2-cycle in $H$ that consists of two multiple edges $e_1$ and $e_2$ joining two black vertices $x_1$ and $x_2$
(see Figure~\ref{f3} (b)),  then  either $int(C)$ or $ext(C)$  contains neither black vertices nor red vertices.
\end{prop}

\proof  As $H$ has no separating 2-cycles,
either $int(C)$ or $ext(C)$ contains no black vertices.
Assume that $int(C)$ does not contain black vertices.

Suppose that $int(C)$ contains red vertices.
Then, $int(C)$ contains white vertices of $D^{\times}$.
As  $int(C)$ does not contain black vertices,
each white vertex in $int(C)$ is of degree at most $2$ in $D^{\times}$.
Thus, we can redraw the edges of $D$ in $int(C)$  such  that these edges make no crossings, and then obtain a 1-planar drawing of $G$ with fewer crossings than $D$, contradicting to the choice of $D$.
Hence $int(C)$ does not contain red vertices
and the conclusion holds.
\proofend

\begin{prop}\label{multiplicity}
{\it Assume that $H$ contains no separating 2-cycles.
Then the edge multiplicity of $H$ is at most 2.}
\end{prop}

\begin{figure}[!ht]
\begin{center}
\includegraphics[width=10cm]{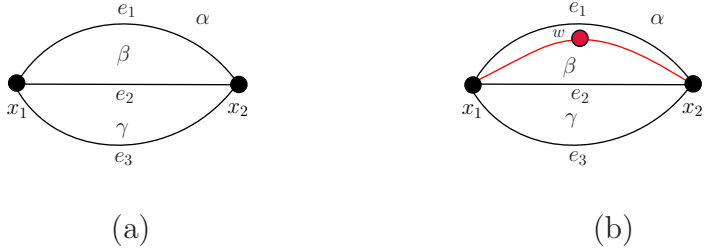}   %{3-1-1}
\vskip 0.2cm

(a) \hspace{5.6 cm} (b)

\caption{Possible  three  multiple edges.
\label{f2}}
\end{center}
\end{figure}

\proof Assume to contrary that $H$ has three  multiple edges $e_1$, $e_2$ and $e_3$ which join the same pair of black vertices $x_1$ and $x_2$.
Then these three edges divide the plane into three regions,  denoted by  $\alpha$, $\beta$ and $\gamma$, as shown in Figure~\ref{f2} (a).
By Proposition \ref{2-cycle},
at least two of these three regions contain neither
red vertices nor black vertices, except on its boundary.
We may assume $\alpha$ and $\gamma$ are such two regions.

Let $P=x_1wx_2$ be the 2-path of $F$
that corresponds to edge $e_3$, where $w$ is a red vertex.
Thus, this path must be within region $\beta$,
as shown in Figure~\ref{f2} (b).

As $P$ is within region $\beta$,
black edges $e_1$ and $e_2$ are in different sets
$int(e_3\cup P)$ and $ext(e_3\cup P)$,
a contradiction to Proposition \ref{empty}.
The proof is then completed.
\proofend

\vskip 0.2cm

An edge of $H$ is called a {\it simple} edge if it is not parallel
to another edge in $H$ and a  {\it  partnered} edge otherwise.
It follows from Proposition \ref{multiplicity} that,
if $H$ has no separating 2-cycles,
then each partnered edge $e$ in $H$ is parallel
to a unique partnered edge $e'$ in $H$.

Let $C$ be a cycle and $P$ be a path in $H$
such that the end vertices of $P$ are the only vertices
in both $C$ and $P$.
When we say that $P$ lies in $int(C)$ (resp. $ext(C)$),
it means that all edges and internal vertices of $P$
lie in $int(C)$ (resp. $ext(C)$).

%\vskip 0.2cm

\noindent
\begin{prop}\label{two 2-paths}
Assume that $H$ has no separating $2$-cycles.
Let $C$ be a $3$-cycle of $H$ consisting of black vertices $x_1$, $x_2$ and $x_3$, and $e$ be the edge on $C$ joining $x_1$ and $x_3$.
Assume that $e'$ is a partnered edge in $H$ which is parallel to $e$.
If $P=x_1wx_3$ and $P'=x_1w'x_3$ are the $2$-paths in $F$
corresponding to $e$ and $e'$ respectively,
then one of $P$ and $P'$ lies in  $int(C)$  and the other in $ext(C)$.
 \end{prop}

\begin{figure}[!ht]
\begin{center}
\includegraphics[width=12cm]{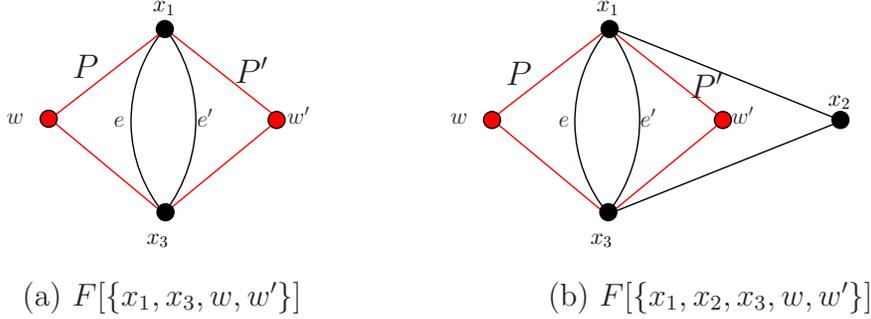} %{8-1}

(a) $F[\{x_1,x_3,w,w'\}]$ \hspace{3 cm}
(b) $F[\{x_1,x_2,x_3,w,w'\}]$

\caption{$x_2$ lies in $ext(C_2)$, where $C_2$ is the cycle $x_1wx_3w'x_1$
%Possible partnered edges  in $H$.
\label{f4}
}
\end{center}
\end{figure}
\vskip 0.2cm

\proof Let $C_0$ denote the 2-cycle of $H$ consisting of edges $e$ and $e'$.
By Proposition \ref{2-cycle},
we may assume that $int(C_0)$ contains neither black vertices nor red vertices.
Thus, both $w$ and $w'$ are in $ext(C_0)$.

Let $C_1$ denote the $3$-cycle of $F$ consisting of edge $e$
and path $P$ and $C'_1$ the $3$-cycle of $F$ consisting of edge $e'$
and path $P'$.
By Proposition \ref{empty},
both $int(C_1)$ and $int(C'_1)$ are empty.
Thus, the subgraph $F[\{x_1,x_3,w,w'\}]$ is as shown in
Figure~\ref{f4} (a).

As these three regions $int(C_0)$, $int(C_1)$ and $int(C'_1)$
do not contain black vertices,
$x_2$ must be in $ext(C_2)$,
where $C_2$ is the $4$-cycle of $F$ consisting of
paths $P=x_1wx_3$ and $P'=x_1w'x_3$.
As $F$ is a plane graph, path $x_1x_2x_3$ must lies in $ext(C_2)$,
as shown in Figure~\ref{f4} (b).

Hence the conclusion holds.
\proofend

\iffalse
\noindent
{\bf Proof}.  By contradiction, there  are two cases:  (1) both  $P$ and $P'$ lie in $ext(C)$;  (2) both  $P$ and $P'$ lie in  $int(C)$. Denote by $P^{*}$ the path $x_1x_2x_3$. First consider the case (1). We first draw $P$ in $ext(C)$ (see the left of Figure 4).  Since  $int(e_1\cup P)$ is ``empty"  by   Proposition \ref{empty},  if $e'_1$ lies in $ext(C)$ then $e'_1$ must  lie in $ext(P\cup P^{*})$. However, at this time no matter how we draw $e'$,  we see that  both
$int(e_1\cup e'_1)$ and $ext(e_1\cup e'_1)$ contain  the red vertex $w$ or  the black vertex $x_2$, a contradiction to Proposition \ref{2-cycle}. Therefore, $e'$ must lie in  $int(C)$ (see the left of Figure 4).
Because $P'$ lies in $ext(C)$  by the assumption in this case,  no matter how $P'$ is drawn, it is seen that $int(e'_1\cup P')$  is not ``empty",   a contradiction to Proposition \ref{empty}.
 Now consider the case (2). Since $P$ lies in $int (C)$ by the assumption in this case, we first draw $P$ in $int (C)$.  Since $int(e_1\cup P)$ is
 ``empty" by  Proposition \ref{empty}, and $P'$ also  lies in $int (C)$ by the assumption,   it  must happen that  $P'$  lies in $int(P\cup P^{*})$ (see the right of Figure 4).   At this time,  no matter how we draw $e'$,  we know that either $int(e'\cup P')$ is not ``empty",  contradicting to Proposition \ref{empty}, or both $int (e_1\cup e'_1)$ and $ext(e_1\cup e'_1)$ contain the  red vertex $w$ or the  black  vertex $x_2$, contradicting to    Proposition \ref{2-cycle}.  This proves the claim.  \hfill{$\square$}

\fi

\begin{prop}\label{simple edge}
Suppose that $H$ has no separating 2-cycles.
For any $3$-cycle $C$ in  $H$,
% consisting of edges  $e_1=x_1x_2$, $e_2=x_2x_3$ and $e_3=x_3x_1$.
if $int(C)$ contains exactly $r$  red vertices,  where $0\le r\le 2$,
then $C$ contains at least $3-r$ simple edges of $F$.
\iffalse
Then  the following hold:
\begin{enumerate}
\item [(1)] if  $int(C)$ contains no   red vertices,  then all the three  edges  $e_i$ $(i=1, 2, 3)$ are  simple;

\item [(3)] if $int(C)$ contains exactly  one  red vertex, then at least two of  the three  edges  $e_i$ $(i=1, 2, 3)$ are  simple;

\item [(3)] if  $int(C)$ contains exactly  two  red vertex, then at least one of  the three  edges  $e_i$ $(i=1, 2, 3)$ are  simple.
\end{enumerate}
\fi
\end{prop}

\proof Let $e_1,e_2$ and $e_3$ be the three edges on $C$.
Suppose that $e_i$ is not a simple edge of $H$, where $1\le i\le 3$.
Then $e_i$ is parallel to another partnered edge $e'_i$ of $H$.
Let $P_i$ and $P'_i$ be 2-paths in $F$ which
correspond to edges $e_i$ and $e'_i$ respectively.
Since $H$ has no separating 2-cycles,
by Proposition \ref{two 2-paths},
$int(C)$ contains a red vertex that is on $P_i$ or $P'_i$.

The above conclusion implies that the number of red vertices in
$int(C)$
is not less than the number of partnered edges on $C$.
Thus, the result holds.
\proofend

\section{Some lemmas \label{sec4}}

Let $G$ be a bipartite graph with partite sets $X$ and $Y$
and ${\cal O}$ be a disk on the plane.
If $D$ is a $1$-planar drawing of $G$
that draws all vertices of $X$ on the boundary of ${\cal O}$ and all vertices of $Y$ and all edges of $G$
in the interior of ${\cal O}$,
%such that each edge is crossed at most once,
then we say that $D$ is a {\it $1$-disc ${\cal O}_X$  drawing} of $G$.
%Clearly, a {\it 1-disc ${\cal O}_X$  drawing} of $G$ is also a 1-planar drawing of  $G$.

\vskip 0.2cm

\begin{lem}\label{disc}
Let $G$ be a bipartite graph with partite sets $X$ and $Y$, and let $D$  be  a 1-disc ${\cal O}_X$ drawing  of $G$  with the minimum number of crossings  $k$. If $|X|=3$, then  $k\in \{0,1,3\}$ and
$|E(G)|\leq 2|Y|+1+\big\lceil \sqrt k\big\rceil$, i.e.,
$$
|E(G)|\leq \ARRAY
{
2|Y|+1, \quad \mbox{if }k=0; \\
2|Y|+2, \quad \mbox{if }k=1; \\
2|Y|+3, \quad \mbox{if }k=3.
}
$$
\end{lem}

\begin{figure}[!ht]
\begin{center}
\includegraphics[width=10cm]{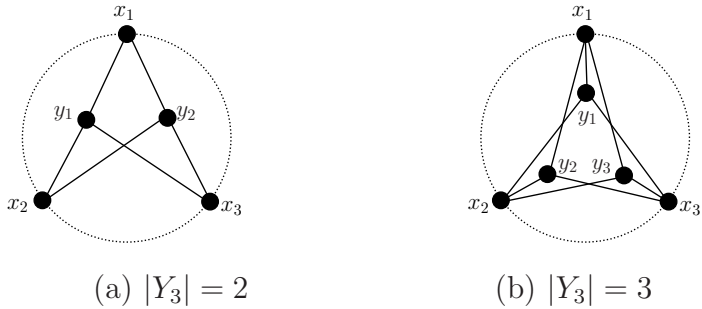}%{9-1}

(a) $|Y_3|=2$ \hspace{3 cm} (b) $|Y_3|=3$

\caption{The 1-disc ${\cal O}_X$ drawing  of $G$
for $Y=Y_3$ and $|Y_3|\in \{2,3\}$
\label{f5}}
\end{center}
\end{figure}

\proof Assume that $|X|=3$.
For any integer $i\ge 0$,
let $Y_i$ be the set of vertices $y$ in $Y$ with $d_G(y)=i$.
As $|X|=3$ and $Y$ is independent in $G$,
$Y_i=\emptyset$ holds for all $i\ge 4$.

It can be checked easily that, for
each vertex $y$ in $\in Y$, if $y\notin Y_3$,
then $y$ is not incident with any crossed edge. Thus,
$G-\bigcup_{i\le 2}Y_i$ has exactly $k$ crossings,
and it suffices to show that  $|Y_3|\le 3$ and
$$
k=\ARRAY
{
0, \quad &\mbox{if } |Y_3|\le 1;\\
1,  &\mbox{if } |Y_3|=2;\\
3,  &\mbox{if } |Y_3|=3.
}
$$

The rest of the proof will be completed by showing the following claims.

\newcommand \aclaim[1]
{\noindent {\bf Claim (#1)}: }

\aclaim{a}  $|Y_3|\le 3$.

Suppose that $|Y_3|\ge 4$.
Then, there exists a bipartite $1$-planar drawing $D'$
isomorphic to $K_{3,2|Y_3|}$
obtained from $D[X\cup Y_3]$
by copying all vertices and edges in the interior of
${\cal O}_X$ to its exterior,
implying that $K_{3,8}$ is 1-planar.
It is a contradiction to the fact that
$K_{3, 7}$ is not  1-planar due to Czap and Hud\'ak~\cite{J.D}.

Thus, Claim (a) holds.

\aclaim{b} If $|Y_3|\le 1$, then $G-\bigcup_{i\le 2}Y_i$ has no
crossings, i.e., $k=0$.

Claim (b) can be verified easily.

\aclaim{c} For any two vertices $y_1,y_2\in Y_3$,
some edge incident with $y_1$
crosses with some edge incident with $y_2$, as shown in
Figure~\ref{f5} (a).

If Claim (c) fails, then $G[X\cup \{y_1,y_2\}]$ is a plane graph
and we can get a drawing of $K_{3,3}$ from $G[X\cup \{y_1,y_2\}]$
by adding a new vertex $y'$ and
three edges joining $y'$ to
all vertices in $X$ in the exterior of ${\cal O}_X$
without any crossing,
implying that $K_{3,3}$ is planar,
a contradiction.

%By Claim 3, the next claim holds.

\aclaim{d} $k= 1$ when $|Y_3|=2$, and $k= 3$ when $|Y_3|=3$.

By Claim (c), $k\ge {|Y_3|\choose 2}$.
By the drawings in Figure~\ref{f5},
$k\le 1$ when $|Y_3|\le 2$, and $k\le 3$ when $|Y_3|\le 3$.
Thus, Claim (d) holds.

The result follows from Claims (a), (b) and (d).
\proofend

\begin{lem}\label{trueface}
Let $G$ be a  plane simple graph with $|V(G)|\geq 3$.
If $G$ has exactly $c$ components and
$t$ $(\geq 0)$ cellular $3$-faces,
then $|E(G)|\leq 2|V(G)|-3-c+\frac{t}{2}$.
\end{lem}

\proof
If $c=1$, since $G$ is simple,
each face of $H$ is a cellular face and has size at least 3.
Then, in this case, the conclusion can be proved easily by applying
the Euler's formula.

Now we assume that $c\ge 2$.
We can obtain a simple and connected plane graph $G'$ from $G$ by
adding $c-1$ edges.

For every noncellular face $F$ of $G$,
we  assume  that its boundary consists of $\ell $  disjoint closed walks of $G$,  and then  we can add $\ell-1$ new edges  (not add  the vertex) by  appropriately drawing these new edges within $F$  so that $F$ is transformed into a cellular face of size at least 4 because $|V(G)|\geq 3$.
Therefore, the resulting  graph $G'$ is a simple and connected plane, and all faces of $G'$ are cellular.

Note that adding the $c-1$ new edges does not produce new cellular 3-faces, and thus $G'$ has exactly $t$ faces of size 3.
The conclusion for connected plane graphs implies that
$$
|E(G')|\leq 2|V(G')|-4+t/2.
$$
As $V(G')=V(G)$ and $|E(G')|=|E(G)|+c-1$, the above inequality
implies that
 $|E(G)|\le 2|V(G)|+t/2-3-c$.
\proofend

\begin{lem}\label{bipartite trueface}
Let $G$ be a simple and bipartite plane graph with $|V(G)|\geq 3$.
If $G$ has exactly $c$ components and
$t$ cellular faces
whose boundaries are of length at least $6$,
then $|E(G)|\leq 2|V(G)|-3-c-t$.
\end{lem}

\proof
If $G$ is connected (i.e. $c=1$), since $G$ is bipartite and simple, then each face of $G$ is  a cellular face, and  has the size at least 4.
Because $G$ has $t$ faces of size at least 6,
it follows from the Euler's formula that
$|E(G)| \leq 2|V(G)|-4-t$.

Now assume that $c\ge 2$.
We can obtain a simple and connected bipartite plane graph
$G'$ from $G$ by adding $c-1$ edges.

For every noncellular face $F$ of $G$ consisting of $\ell$
distinct closed walks,
similar to the proof of Lemma \ref{trueface},
we can add $\ell-1$ new edges within this
noncellular face so that $F$ is transformed into
a cellular face.
We can ensure that those new added edges join the vertices  in  different partite sets of $G$.
Hence the resulting plane graph $G'$ is simple, bipartite
and connected.
Clearly, all faces of $G'$ are cellular,
and $G'$ has at least $t$ faces whose boundaries
are of length at least $6$.
The conclusion for bipartite and connected plane graphs implies that
$$
|E(G')|\leq 2|V(G')|-4-t.
$$
As $V(G')=V(G)$ and $|E(G')|=|E(G)|+c-1$, the above inequality
implies that $|E(G)|\le 2|V(G)|-3-c-t$.
\proofend

\noindent {\bf Remark}:
Lemmas~\ref{trueface} and~\ref{bipartite trueface}
can be strengthened when $G$ contains isolated vertices.
Let $V_{\ge 1}(G)$ be the set of non-isolated vertices in $G$.
Then, under the condition  $|V_{\ge 1}(G)|\ge 3$,
the conclusions of both Lemmas~\ref{trueface} and~\ref{bipartite trueface}
still hold after $|V(G)|$ is replaced by $|V_{\ge 1}(G)|$.

\iffalse
\vskip 0.2cm
\noindent
{\bf Proof}.  If $G$ is connected, since $G$ is bipartite and simple, then each face of $G$ is  a cellular face, and  has the size at least 4.  Because $G$ has $t$ faces of size 6, it follows from the Euler's formula that
$|E(G)| \leq 2|V(G)|-4-t$. If $G$ is not connected, for every noncellular face, we can do the analogy  like as  in the proof of Lemma \ref{trueface}, by adding a number of new edges within every noncellular faces of $G$ so as to change this face into a  new cellular face.  We can ensure that those new added edges join the vertices  in  different partite sets of $G$.  Hence  the resulting  plane graph $G'$ is still a  bipartite, connected and simple graph.  Clearly, all faces of $G'$ are cellular,  and $G'$  has at least  $t$ faces of sizes 6.     Therefore, by the Euler's formula $|E(G)|\leq |E(G')|\leq 2|V(G)|-4-t$.
\hfill{$\square$}

\fi

\section{Proof of Theorem~\ref{th-main}
\label{sec5}}

The whole section contributes to the proof of Theorem~\ref{th-main}.

\vspace{0.2 cm}

\noindent {\it Proof} of Theorem~\ref{th-main}.
Suppose that Theorem~\ref{th-main} fails and
$\chi$ is the minimum integer with $\chi\ge 2$
such that for some bipartite 1-planar
graph $G$ with partite sets $X$ and $Y$,
where $\chi=|X|\le |Y|$,
$|E(G)|>2|V(G)|+4|X|-12$ holds.

We will prove the following claims to show that this assumption
leads to a contradiction.

\inclaim $\chi\ge 4$.
\def \cli {1}

\proofn
Let $G$ be any bipartite $1$-planar graph with bipartitions $X$ and $Y$,
where $2\le |X|\le |Y|$.
If $|X|=2$, obviously,
$$
|E(G)|\le 2|Y|=2(2+|Y|)+4\times 2-12
=2|V(G)|+4|X|-12.
$$
Now assume that $|X|=3$.
Let $Y_i$ be the set of vertices $y$ in $Y$ with $d_G(y)=i$.
Then $|E(G)|\le  2|Y|+|Y_3|$.
Since the complete bipartite graph $K_{3,7}$
is not 1-planar (see \cite{J.D}),
we have $|Y_3|\le 6$, implying that $|E(G)|\le  2|Y|+6$.
As $x=3$, we have
$$
2|V(G)|+4|X|-12=6|X|+2|Y|-12=2|Y|+6.
$$
Thus, $|E(G)|\le 2|V(G)|+4|X|-12$.

By the assumption of $\chi$, we have $\chi\ge 4$.
\proofend

\iffalse
Assume that $x=|X|\geq 4$. By  induction hypothesis, for  a bipartite 1-planar graph $G'$ (connected or disconnected) with partite sets $X'$ and $Y'$ such that $2\leq |X'|\leq |Y'|$, if $|X'| < x=|X|$, then $|E(G')|\leq 2|V(G')|+4|X'|-12$.
\fi

In the following, we assume that $G$ is a bipartite $1$-planar graph
with bipartitions $X$ and $Y$, where $\chi=|X|\le |Y|$,
such that
\begin{equation}\label{eq0}
|E(G)|> 2|V(G)|+4|X|-12.
\end{equation}

Let $D$ be a 1-planar drawing of $G$ with the minimum number of crossings
and $W$ be the set of its crossings.
Introduced in Section~\ref{sec3},
$D^{\times}_W$ is a plane graph extended from $D^{\times}$,
and $F$ and $H$ are the subgraphs $D^{\times}_W[X\cup W]$
and $D^{\times}_W[X]$ of $D^{\times}_W$ respectively.
All vertices in $X$ are black vertices,
all vertices in $Y$ are white vertices
and all vertices in $W$ are red vertices.

We are now going to prove the following claim.

\inclaim $H$ has no separating 2-cycles.
\def \clii {2}

\proofn Assume to the contrary that  $H$ has a separating 2-cycle $C$
consisting of two parallel edges $e_1$ and $e_2$
joining black vertices $x_1$ and $x_2$
(see Figure~\ref{f3} (b)),
such that both $int(C)$ and  $ext(C)$ contain black vertices.

Let $G_1=G\bigcap INT(C)$ and $G_2=G\bigcap EXT(C)$.
Obviously, both $G_1$ and $G_2$ are  bipartite 1-planar  subgraphs of $G$.  Moreover, we can see that $G_1\cup G_2= G$, and $ G_1\cap G_2=\{x_1, x_2\}$.

For $i=1,2$, $G_i$ has a bipartition $X_i$ and $Y_i$,
where $X_i=X\cap V(G_i)$ and $Y_i=Y\cap V(G_i)$.
Clearly, $|X_1|+|X_2|=|X|+2=\chi+2$ and $|Y_1|+|Y_2|=|Y|$.

Since $C$ is a separating cycle of $H$,
both $int(C)$ and $ext(C)$  contain black vertices,
implying that $|X_i|\ge 3$ for both $i=1,2$.
As $|X_1|+|X_2|=|X|+2$, we have $|X_i|<|X|=\chi$
and so $\min\{|X_i|,|Y_i|\}\le |X_i|<\chi$.

For $i=1,2$, if $|Y_i|\le 1$, then
$|E(G_i)|\le |Y_i|\cdot |X_i|$ and
$|E(G_i)|\le 2|V(G_i)|+4|X_i|-12$ holds;
if $|Y_i|\ge 2$, then $2\le \min\{|X_i|,|Y_i|\}\le |X_i|<\chi$
and the assumption on $\chi$ implies that
the conclusion of Theorem~\ref{th-main} holds for $G_i$,
i.e.,
\begin{eqnarray}\label{eq10}
|E(G_i)|\le 2|V(G_i)|+4\min\{|X_i|,|Y_i|\}-12
\le 2|V(G_i)|+4|X_i|-12.
\end{eqnarray}
Thus, by (\ref{eq10}),
\begin{eqnarray}\label{eq1}
|E(G)|&=&|E(G_1)|+|E(G_2)|
\nonumber \\
&\le &2(|V(G_1)|+|V(G_2)|)+4(|X_1|+|X_2|)-24
\nonumber \\
&=&2(|V(G)|+2)+4(|X|+2)-24
\nonumber \\
&=&2|V(G)|+4|X|-12,
\end{eqnarray}
which contradicts to the assumption in (\ref{eq0}).

Hence Claim~\clii\ holds.
\proofend

It is known from  Proposition \ref{multiplicity} that
the edge multiplicity of each edge in $H$ is  at most 2.
Then, there exists a subset $A$ of $E(H)$ such that
both $H\langle A\rangle $ and $H-A$ are simple graphs and
each edge in $A$ is parallel to some edge in $E(H)-A$,
where $H\langle A\rangle $ is the spanning subgraph of $H$
with edge set $A$ and $H-A$ is the graph obtained from $H$
by removing all edges in $A$.
Clearly, $|A|$ is the number of pairs of edges $e$ and $e'$ in $H$
which are parallel.

Let $H'$ denote $H-A$.
Obviously, $|A|\leq |E(H')|$, and
  \begin{equation}\label{partened edges}
   |E(H)|=|E(H')|+|A|.
\end{equation}

\inclaim $H'$ contains at least one cellular 3-face.
\def \cliii {3}

\proofn Suppose that $H'$ has no cellular 3-faces.
As $H'$ is a simple plane graph,
by Lemma \ref{trueface}, $|E(H')|\leq 2|V(H')|-4=2|X|-4$.
Because  $|A|\leq |E(H')|$, it follows from \eqref{partened edges} that $|E(H)|\leq 2|E(H')|\leq 4|X|-8$.

Since  each edge of $H$ is in 1-1 correspondence with  a crossing of
a drawing $D$ of $G$,
we can obtain a simple bipartite plane graph (possibly disconnected),
denoted by $G'$,
by removing $|E(H)|$ edges from $G$ each of which
is a crossed edge of $G$.
By Lemma~\ref{bipartite trueface},
$|E(G')|\leq 2|V(G')|-4= 2|V(G)|-4$.
Therefore,
$$
|E(G)|=|E(G')|+|E(H)|
\le |E(G')|+4|X|-8=2|V(G)|+4|X|-12,
$$
a contradiction to the assumption in (\ref{eq0}).

Hence Claim~\cliii\ holds.
\proofend

Now we assume that $H'$ has  exactly  $t$ cellular 3-faces,
where $t\geq 1$.
Let $\mathscr{T}(H')$ denote the set of  cellular 3-faces in $H'$.
So $t=|\mathscr{T}(H')|$.

For each  $\Delta\in \mathscr{T}(H')$,
for convenience we also use ``$\Delta$'' to represent the 3-cycle corresponding the boundary of $\Delta$
if there is no confusions in the context.
Let $G_{\Delta}=G\bigcap INT(\Delta)$.
Since $\Delta$ is a cellular 3-face of  $H'$,
there are no black vertices lying in $int(\Delta)$,
and thus $G_{\Delta}$ is a bipartite with
with exactly three black vertices,
which lie on the boundary of the face $\Delta$.

Let $\Delta\in \mathscr{T}(H')$.
Since $D$ is a 1-planar drawing of $G$ with minimal number of crossings,
the induced subdrawing of $D$ of $G_{\Delta}$
is a 1-disc ${\cal O}_{X_{\Delta}}$ drawing
of $G_{\Delta}$ with the minimum number of crossings,
where $X_{\Delta}$ is the set of three black vertices
on the boundary of $\Delta$.
Otherwise, we redraw the edges of $G$  lying in the interior of $\Delta$, and obtain a 1-planar drawing of $G$ with fewer crossings than $D$, contradicting to the choice of $D$.
By Lemma \ref{disc},  the number of crossings of $D$
in $int (\Delta)$ is a number in the set $\{0, 1,3\}$.

For any $j\in \{0,1,3\}$,
let $\mathscr{T}^{(j)}(H')$
be the set of members $\Delta$ in $\mathscr{T}(H')$
such that $int (\Delta)$ contains exactly $j$ crossings of $D$.
Assume that
$\mathscr{T}^{(j)}(H')=\{\Delta^{(j)}_i: 1\le i\le t_j\}$,
where $t_j=|\mathscr{T}^{(j)}(H')|$.
Thus, $t_0+t_1+t_3=t$.

For each $\Delta^{(j)}_i\in \mathscr{T}^{(j)}(H')$,
let $e^{(j)}_i$ be the number of the edges of
the graph $G_{\Delta^{(j)}_i}$
and $y^{(j)}_i$
be the number of white vertices  in $int(\Delta^{(j)}_i)$.

\inclaim $\sum\limits_{j\in \{0,1,3\}}
\sum\limits_{i=1}^{t_{j}}e^{(j)}_i
\leq 2\sum\limits_{j\in \{0,1,3\}}
\sum\limits_{i=1}^{t_{j}}y^{(j)}_i +(t_0+2t_1+3t_3)$.

\def \cliv {4}

\proofn
By Lemma~\ref{disc},
$e^{(j)}_i\leq 2y^{(j)}_i+1+\big\lceil \sqrt {j}\big\rceil$ holds
for any $j\in \{0,1,3\}$ and $1\le i\le t_j$.
Thus,
\iffalse
\begin{equation}\label{interior edges}
\sum\limits_{j\in \{0,1,3\}}
\sum\limits_{i=1}^{t_{j}}e^{(j)}_i
\leq 2\sum\limits_{j\in \{0,1,3\}}
\sum\limits_{i=1}^{t_{j}}y^{(j)}_i +(t_0+2t_1+3t_3).
\end{equation}
\fi
Claim~\cliv\ holds.
\proofend

\inclaim
$|E(H)|\leq 4|X|-8+t-({3t_0}+2t_1)/2$.

\def \clv {5}

\proofn
Since $H'$ is a simple plane graph with exactly $t\geq 1$ cellular 3-faces
and $|V(H')|=|X|=\chi\ge 4$,
by Lemma \ref{trueface},
\begin{equation}\label{H' edges}
|E(H')|\leq 2|V(H')|-4+\frac{t}{2}=2|X|-4+\frac{t}{2}.
\end{equation}

On the other hand,  by Proposition \ref{simple edge},
for any $j\in \{0,1,3\}$ and $1\le i\le t_j$,
at least $3-j$ edges
on the boundary of $\Delta^{(j)}_i$  are simple edges,
implying that at least $3-j$ edges
on the boundary of $\Delta^{(j)}_i$ are in $H'$.
Because each simple edge of $H'$ belongs to
the boundaries of  at most two different faces of $H'$,
it follows that $|E(H')|\ge (3t_{0}+2t_1)/2$.
Then, by (\ref{partened edges}),
$$
|A|\leq |E(H')|-(3t_{0}+2t_1)/2,
$$
and therefore, by \eqref{partened edges} and \eqref{H' edges},
\begin{equation}\label{H edge}
|E(H)|=|E(H')|+|A| \leq 4|X|-8+t-(3t_{0}+2t_1)/2.
\end{equation}
Thus, Claim~\clv\ holds.
\proofend

Let $D'$ denote the drawing obtained from $D$ by deleting
all white vertices and edges of $D$ that
lie in the interiors of all cellular 3-faces $\Delta^{(j)}_i$ of $H'$,
where $j\in \{0,1,3\}$ and $1\leq i\leq t_j$,
and let $G'$ denote the graph represented by $D'$.

We see that the graph $G'$ is a bipartite 1-planar graph
with a bipartition $X$ and $Y'=Y\cap V(G')$,
where $|Y'|=|Y|-\sum\limits_{j\in \{0,1,3\}}
\sum\limits_{i=1}^{t_{j}}y^{(j)}_i $. Thus,
 \begin{equation}\label{G' vertex}
|V(G')|=|V(G)|-\sum_{j\in \{0,1,3\}}\sum\limits_{i=1}^{t_{j}}y^{(j)}_i.
\end{equation}
 and
\begin{equation}\label{G' edge}
|E(G')|=|E(G)|-\sum_{j\in \{0,1,3\}}\sum\limits_{i=1}^{t_{j}}e^{(j)}_i.
\end{equation}
As the number of crossings of $D$ equals to $|E(H)|$
and $D'$ has no crossings lying  in the interior of any cellular 3-face of $H'$, $D'$ has exactly  $|E(H)|-(t_1+3t_3)$ crossings.
\iffalse
, moreover that $G'$
has one partite set with $x$ black vertices, the other partite set with  $\Big(y-\sum\limits_{j=0}^{j=3}\sum\limits_{i=1}^{t_{j}}y^{(j)}_i\Big)$ white vertices. Note that a red vertex is in 1-1 correspondence with a crossing of $D$ and the number of crossings of $D$ equals to $|E(H)|$ (also, equals to the  number of red vertices).
Because any cellular 3-face of $H'$ contains no edges of $G'$, we know that  $D'$ has no crossings  lying  in the interior of a cellular  3-face of $H'$. Therefore,  $D'$ has in all  $\Big(|E(H)|-(0+t_1+2t_2+3t_3)\Big)$ crossings.  Obviously,
\fi

For each crossing of $D'$, we remove exactly one crossed edge from $G'$
and obtain a bipartite plane graph $G^{*}$.
Thus, $|E(G^*)|=|E(G')|-(|E(H)|-(t_1+3t_3))$.
Then, \eqref{G' edge} implies that
\begin{eqnarray}\label{edge-G}
|E(G^{*})|
&=&\Big(|E(G)|-\sum_{j\in \{0,1,3\}}\sum\limits_{i=1}^{t_{j}}e^{(j)}_i\Big)-|E(H)|+(t_1+3t_3).
\end{eqnarray}
Clearly, by \eqref{G' vertex},
\begin{equation}\label{G* vertex}
|V(G)|=|V(G')|+\sum_{j\in \{0,1,3\}}\sum\limits_{i=1}^{t_{j}}y^{(j)}_i
=|V(G^*)|+\sum_{j\in \{0,1,3\}}\sum\limits_{i=1}^{t_{j}}y^{(j)}_i.
\end{equation}

Now,  we shall obtain an upper bound of $|E(G^{*})|$
in terms of $|V(G^*)|$
by  constructing a bipartite plane graph with
at least $t$ cellular  6-faces.

\inclaim $|E(G^{*})|\leq  2|V(G^{*})|-4-t$.
\def \clvi {6}

\proofn Note that the simple and bipartite plane graph $G^{*}$
is obtained from $G$ by removing all white vertices and edges of $G$
lying in the interiors of all cellular 3-faces of $H'$
and, for each crossing of $D$ not lying in
any cellular 3-face of $H'$,
removing exactly one edge of $G$ involved in this crossing.

Now let $G^{**}$ denote the graph obtained from $G^*$
by adding all black edges in $H'$ which belong to the boundary
of cellular 3-faces of $H'$
and then subdividing each of these added edges.
Let $m$ be the number of edges in $H'$
that belong to the boundaries of cellular 3-faces of $H'$.
Then
\begin{equation}\label{G**}
|V(G^{**})|=|V(G^{*})|+m \quad \mbox{and}\quad
|E(G^{**})|=|E(G^{*})|+2m.
\end{equation}

\iffalse the following two  procedures:
\begin{enumerate}
\item[(1)] keep all the edges belonging to the boundary of
 a cellular 3-face of $H'$,  and  retain all of the cellular  3-faces $H'$,  intuitively, at this  moment the resulting plane graph  is: $G^{*}\bigcup \Big(\sum\limits_{j=0}^{j=3}\sum\limits_{i=1}^{t_{j}}\Delta^{(j)}_i\Big)$;
\item[(2)] for each edge belonging to the boundary of a  cellular 3-face of $H'$,  subdivide this edge by inserting a new white vertex.
\end{enumerate}
\fi

Because the edges of $H$ (and thus $H'$) are not crossed with the edges of $G$ (and thus $G^{*}$),
we observe that $G^{**}$ is also a simple and bipartite plane graph
and has at least $t$ cellular  6-faces.
Applying Lemma \ref{bipartite trueface} to $G^{**}$
yields that
$$
|E(G^{**})| \leq 2|V(G^{**})|-4-t.
$$
Then, (\ref{G**}) implies that
$|E(G^{*})|\leq  2|V(G^{*})|-4-t$.
This proves the claim.
\proofend

\inclaim $|E(G)|\le 2|V(G)|+4|X|-12-t_0/2$.
\def \clvii {7}

\proof By (\ref{edge-G}),
we have
$$
|E(G)|=|E(G^{*})|+|E(H)|
+\sum\limits_{j\in \{0,1,3\}}\sum\limits_{i=1}^{t_{j}}e^{(j)}_i
-(t_1+3t_3).
$$
Then, by Claims~\cliv, \clv~and~\clvi,
\begin{eqnarray*}
|E(G)|&\leq &\Big(2|V(G^{*})|-4-t\Big)+
\Big (4|X|-8+t-(3t_0+2t_1)/2\Big )
\\
& &+\Big(2\sum_{j\in \{0,1,3\}}
\sum_{i=1}^{t_j}
y^{(j)}_i +(t_0+2t_1+3t_3)\Big)
-(t_1+3t_3) \\
&=&2|V(G)|+4|X|-12-t_0/2
\hspace{5 cm} \mbox{by } (\ref{G* vertex})  \\
     &\leq& 2|V(G)|+4|X|-12.
\end{eqnarray*}
{}
\proofend

Clearly, Claim 7 contradicts the assumption in (\ref{eq0}).
Hence Theorem~\ref{th-main} holds.
\proofend

\section{Remarks\label{sec6}}

For any $x\ge 3$ and $y\ge 6x-12$,
Czap, Przybylo  and \u{S}krabul\'{a}kov\'{a}
\cite[Lemma 4]{CPS}
constructed a bipartite $1$-planar graph $G$ with partite sets $X$ and $Y$
 such that $|E(G)|=2|V(G)|+4|X|-12$.
Notice that the $1$-planar drawing $D$ of this graph $G$
given in \cite{CPS}
has the following property:
\begin{quote}
(*) each vertex in $X$ is incident with crossed edges in $D$.
\end{quote}
The proof of Theorem~\ref{th-main} also yields that,
if $|E(G)|=2|V(G)|+4|X|-12$ holds for a
bipartite $1$-planar graph $G$ with  partite sets
$X$ and $Y$, where $4\le |X|\le |Y|$,
and $D$ is a $1$-planar drawing of $G$ with the minimum number of
crossings,
then the graph $H'$ introduced in Section~\ref{sec3} does not have isolated vertices, i.e., property (*) above holds.

Based on the above observations, we propose the following
problem.

\begin{prob}\label{prob1}
For any bipartite $1$-planar graph $G$ with partite sets
$X$ and $Y$, where $4\le |X|\le |Y|$,
if $|E(G)|=2|V(G)|+4|X|-12$, does
property (*) hold for every $1$-planar drawing $D$ of $G$
with the minimum number of crossings?
\end{prob}

From Claims 3 and 5 in the proof of Theorem~\ref{th-main},
we can see that if $H'$ does not have separating $2$-cycles
and $|X_{>0}|\ge 3$,
where $X_{>0}$ is the set of non-isolated vertices in $H'$
(i.e., the set of vertices in $X$
which are incident with crossed edges of $D$),
then $|E(G)|\le 2|V(G)|+4|X_{>0}|-12$ holds.

\begin{prob}\label{prob2}
Let $G$ be a bipartite $1$-planar graph with  partite sets
$X$ and $Y$, where $4\le |X|\le |Y|$.
If $D$ is a $1$-planar drawing of $G$
with the minimum number of crossings and
$|X_{>0}|\ge 3$, where $X_{>0}$ is the set of vertices in $X$
which are incident with crossed edges of $D$,
does $|E(G)|\le 2|V(G)|+4|X_{>0}|-12$ hold?
\end{prob}

Theorem~\ref{th-main} shows that $|E(G)|\le 2|V(G)|+4x-12$ holds
for any bipartite $1$-planar graph $G$ with bipartite sets of sizes
$x$ and $y$, where $2\le x\le y$.
For any $x\ge 3$ and $y\ge 6x-12$,
Czap, Przybylo and \u{S}krabul\'{a}kov\'{a}~\cite{CPS}
constructed a bipartite $1$-planar graph $G$
with bipartite sets of sizes $x$ and $y$ and
$|E(G)|=2|V(G)|+4x-12$.
Notice that these graphs constructed in~\cite{CPS}
have minimum degree $3$.
By Theorem~\ref{th1},  any bipartite $1$-planar graph
of $n$ vertices has at most $3n-8$ edges,
implying that its minimum degree is at most $5$.
We wonder if the result in Theorem~\ref{th-main}
can be improved for bipartite $1$-planar graphs
with higher minimum degrees or connectivity.

\begin{prob}\label{prob3}
Let $4\le t\le 5$ and
$G$ be any bipartite $1$-planar graph with partite sets
$X$ and $Y$, where $t\le |X|\le |Y|$.
If $G$ is $t$-connected (or $\delta(G)=t$),
does $|E(G)|\le 2|V(G)|+f(t)|X|+c$ holds for some $f(t)< 4$?
\end{prob}

Let $t\ge 2$.
A drawing of a graph is $t$-{\it planar} if each of its edges is
crossed at most $t$ times.
If a graph has a $t$-planar drawing, then it is $t$-{\it planar}.
Does Theorem~\ref{th-main} have an analogous result
for bipartite $2$-planar graphs?

\begin{prob}\label{prob4}
Let $G$ be a bipartite $2$-planar graph with partite sets
$X$ and $Y$, where $2\le |X|\le |Y|$.
Determine constants $a,b$ and $c$ such that
$|E(G)|\le a|V(G)|+b|X|+c$.
\end{prob}

Lemma~\ref{disc} gives an upper bound for the size of a bipartite graph
$G$ with partite sets $X$ and $Y$, where $|X|=3$,
which has a 1-disc ${\cal O}_X$ drawing.
Can this result be generalized for such a bipartite graph
without the condition that $|X|=3$?

\begin{prob}\label{prob5}
Let $G$ be a bipartite graph with partite sets $X$ and $Y$
which has a 1-disc ${\cal O}_X$ drawing.
Is it true that $|E(G)|\leq 2|Y|+5|X|/3-2$?
\end{prob}


\begin{thebibliography}{99}\setlength{\itemsep}{-2mm}




\rebibitem{JAB} J. A. Bondy and U. S. R. Murty, {\em Graph Theory},
Grad. Texts in Math., vol.244,
Springer, New York, 2008.

\rebibitem{HR}  H. Bodendiek, R. Schumacher and K. Wagner,
\H{U}ber 1-optimale Graphen, {\em Math. Nachr.}
{\bf 117} (1984), 323-339.

\rebibitem{J.D} J. Czap and D. Hud$\acute{a}$k, 1-planarity of complete multipartite graphs, {\em Discrete Appl. Math.} {\bf 160} (2012), 505-512.

\rebibitem{JDT}J. Czap and D. Hud$\acute{a}$k and T. Madras, Joins of 1-planar graphs, {\em Acta  Mathematica Sinica, English Series}
{\bf 30}(11) (2014), 1867-1876.

\rebibitem{CPS} J. Czap, J. Przybylo,  and E. \u{S}krabul\'{a}kov\'{a}, On an extremal problem in the class of bipartite 1-planar graphs, {\em Discuss. Math., Graph Theory}  {\bf 36} (2016), 141-151.

\rebibitem{IT} I. Fabrici and T. Madars, The structure of 1-planar graphs, {\em Discrete Math.} {\bf 307} (2007), 854-865.

\rebibitem{DVK} D.V. Karpov, Upper bound on the number of edges of an almost planar bipartite graph, {\em J. Math. Sci.}
{\bf 196} (2014), 737-706.

\rebibitem{SK} S. G. Kobourov, G. Liotta and F. Montecchiani, An annotated bibliography on 1-planarity, {\em Comput. Sci. Rev.}
{\bf 25} (2017),  49-67.


\rebibitem{DJK} D.J. Kleitman, The crossing number of $K_{5,n}$,
{\em J.  Combin. Theory} {\bf 9} (1970), 315-323.

\rebibitem{JP} J. Pach and G. T\'{o}th,
Graphs drawn with few crossings per edge, {\em Combinatorica}
{\bf 17}(3) (1997), 427-439.

\rebibitem{GR} G. Ringel, Ein Sechsfarbenproblem auf der Kugel, {\em Abh. Math. Sem. Univ. Hamburg}
{\bf 29} (1965) 107-117.


\rebibitem{S} E. Sopena, personal communication, Seventh Cracow Conference in Graph Theory, Rytro (2014).

\rebibitem{XZ} X. Zhang, Drawing complete multipartite graphs on the plane with restrictions on
crossings, {\em Acta Mathematica Sinica, English Series}
{\bf 30} (12) (2014), 2045-2053.












%\vskip 1cm

%\bibitem{MA} M. O. Albertson, Chromatic number, independence ratio, and crossing number,
%{\em Ars Math. Contemp}. {\bf 1} (1) (2008) 1-6.

%\bibitem{CB} C. Bachmaier, F. J. Brandenburg, K.Hanauer, et al, NIC-planar graphs. {\em Discrete Applied Mathematics}, {\bf 232} (2017) 23-40.

%\bibitem{JB} J. Bar\'{a}t, G T\'{o}th. Improvements on the density of maximal 1-planar graphs. {\em Journal of Graph Theory}, {\bf88}(1) (2018), 101-109.



%\bibitem{TB} T. Biedl, A note on 1-planar graphs with minimum degree 7, arXiv preprint arXiv:1910.01683, 2019.



%\bibitem{FJ} F. J. Brandenburg, D. Eppstein, A. Glei{\ss}ner, et al, On the density of maximal 1-planar graphs. International Symposium on Graph Drawing. Springer, Berlin, Heidelberg, 2012, 327-338.


%\bibitem{JC} J. Czap, D. Hud\'{a}k, On drawings and decompositions of 1-planar graphs. {\em The Electronic Journal of Combinatorics}, (2013) P54-P54.

%\bibitem{WD} W. Didimo, Density of straight-line 1-planar graph drawings, {\em Information Processing Letters}, {\bf 113} (7) (2013) 236-240.



%\bibitem{DH} D. Hud\'{a}k, T. Madaras, Y. Suzuki. On properties of maximal 1-planar graphs, {\em Discuss. Math., Graph Theory}, {\bf 32}, (2012) 737-747.

%\bibitem{HS} H. Schumacher, Zur struktur 1-planarer graphen, {\em Math. Nachr}. {\em 125} (1986) 291-300.

%\bibitem{YS} Y. Suzuki, Re-embeddings of maximum 1-planar graphs, {\em SIAM J. Discrete Math}. {\bf 24} (4) (2010) 1527-1540.



%\bibitem{XZ2} X. Zhang, G. Liu, The structure of plane graphs with independent crossings and
%its application to coloring problems, {\em Central Europ. J. Math} {\bf 11} (2) (2013) 308-321.

\end{thebibliography}
\end{document}